\documentclass[11pt]{article}
\newcommand{\bbR}{{\rm I\hspace{-0.7mm}R}}

\usepackage{latexsym}

\begin{document}

\centerline{\bf \large On A Theorem In Multi-Parameter Potential Theory }

\vskip .1in

\centerline{\tt By MING YANG}

\vskip .1in

\indent{\footnotesize
Let $X$ be an $N$-parameter additive L\'evy process in $\bbR^d$ with
L\'evy exponent $(\Psi_1,\cdots,\Psi_N)$ and
let $\lambda_d$ denote Lebesgue measure in $\bbR^d.$
We show that
$$E\{\lambda_d(X(\bbR^N_+))\}>0\Longleftrightarrow
\int_{\bbR^d}
\prod_{j=1}^N\mbox{Re}\left(\frac{1}{1+\Psi_j(\xi)}\right)
d\xi<\infty.$$ This was previously proved by 
Khoshnevisan, Xiao and Zhong [1] under a sector condition.}

\vskip .1in
\noindent
{\small \noindent 2000 {\it Mathematics Subject Classification}. Primary
60G60, 60G51; secondary 60G17.
  
\noindent
{\it Key words and phrases}. Additive L\'evy processes, Hausdorff dimension, multiple points.}

\vskip .15in
\centerline{\textsc{1. Introduction and Proof}}
\vskip .15in

Let
$X^1_{t_1},~X^2_{t_2},\cdots,X^N_{t_N}$ be $N$ independent L\'evy
processes in $\bbR^d$ with their respective L\'evy exponents 
$\Psi_j,~j=1,2,\cdots,N$. The random field
$$X_t=X^1_{t_1}+X^2_{t_2}+\cdots+X^N_{t_N},~~~~~t=(t_1,t_2,\cdots,t_N)\in\bbR^N_+$$
is called the additive
L\'evy process. Let $\lambda_d$ denote Lebesgue measure in $\bbR^d.$

\vskip 0.15in
\noindent
{\bf Theorem 1.1}
{\it
Let $X$ be an additive L\'evy process in $\bbR^d$ with
L\'evy exponent $(\Psi_1,\cdots,\Psi_N).$
Then}
$$E\{\lambda_d(X(\bbR^N_+))\}>0\Longleftrightarrow
\int_{\bbR^d}
\prod_{j=1}^N\mbox{Re}\left(\frac{1}{1+\Psi_j(\xi)}\right)
d\xi<\infty.\eqno(1.1)$$

\vskip .1in
 
\noindent
Recently,
Khoshnevisan, Xiao and Zhong [1] proved that if 
$$\mbox{Re}\left(\prod_{j=1}^N
\frac{1}{1+\Psi_j(\xi)}\right)\ge\theta\prod_{j=1}^N
\mbox{Re}\left(\frac{1}{1+\Psi_j(\xi)}\right)\eqno(1.2)$$
for some constant $\theta>0$ then
Theorem 1.1 holds. In fact the proof of Theorem 1.1 does not need any condition.

\vskip .15in
\noindent
{\bf Proof of Theorem 1.1:} Define
$$\mathcal{E}_{\Psi}(\mu)=(2\pi)^{-d}\int_{\bbR^d}|\hat{\mu}(\xi)|^2
\prod_{j=1}^N\mbox{Re}\left(\frac{1}{1+\Psi_j(\xi)}\right)
d\xi$$
where $\mu$ is a probability measure on a compact set $F\subset\bbR^d$
and $\hat{\mu}(\xi)=\int_{\bbR^d}e^{i\xi\cdot x}\mu(dx).$ Let $F=\{0\}\subset\bbR^d$ and
$\delta_0$ be the point mass at $0\in\bbR^d.$ We first quote a key lemma of [1]:

\vskip .15in
\noindent
{\bf Lemma 5.5}~~
{\it Suppose $X$ is an additive L\'evy process in $\bbR^d$ that satisfies Condition (1.3),
and that $\int_{\bbR^d}\prod_{j=1}^N|1+\Psi_j(\xi)|^{-1}d\xi<+\infty,$
where $\Psi=(\Psi_1,\cdots,\Psi_N)$ denotes the L\'evy exponent of $X$.
Then, for all compact sets $F\subset\bbR^d$, and for all $r>0,$}
$$E\{\lambda_d(X([0,r]^N\oplus F)\}\le \theta^{-2}(4e^{2r})^N\cdot\mathcal{C}_{\Psi}(F),$$ 
{\it where $\theta>0$ is the constant in Condition (1.3).}

\vskip .15in

By reviewing the whole process of the proof of Theorem 1.1 of [1] given by  
Khoshnevisan, Xiao and Zhong, our Theorem 1.1 certainly follows if we instead prove
the following statement:

\vskip .15in

{\it Let $X$ be any additive L\'evy process in $\bbR^d.$ If
$\int_{\bbR^d}\prod_{j=1}^N|1+\Psi_j(\xi)|^{-1}d\xi<+\infty,$
then}
$$E\{\lambda_d(X([0,r]^N))\}\le\frac{c_{N,d,r}}{\mathcal{E}_{\Psi}(\delta_0)}\eqno(1.3)$$
{\it for some constant $c_{N,d,r}\in(0,\infty)$ depending on $N,~d,~r$ only.}

\vskip .15in
\noindent
Clearly, all we have to do is to complete Eq. (5.11) of [1] without 
bothering ourselves with Condition (1.3) of [1].
Since $\delta_0$ is the only probability measure on 
$F=\{0\}$, letting $\eta\rightarrow0,~k\rightarrow\infty,$ and $\varepsilon\rightarrow0$ 
and using the integrability condition
$\int_{\bbR^d}\prod_{j=1}^N|1+\Psi_j(\xi)|^{-1}d\xi<+\infty$
yield
$$\mathcal{E}_{\Psi}(\delta_0)\ge c_1\left|\int_{\bbR^d}\mbox{Re}\left(\prod_{i=1}^N
\frac{1}{1+\Psi_i(\xi)}\right)d\xi\right|^2E\{\lambda_d(X([0,r]^N))\}\eqno(1.4)$$
where $c_1\in(0,\infty)$ is a constant depending on $N,~d,~r$ only.

Consider the $2^{N-1}$ similar additive L\'evy processes (including $X_t$ itself)
$X^{\pm}_t=X^1_{t_1}\pm X^2_{t_2}\pm\cdots\pm X^N_{t_N}.$
Here, $\pm$ is merely a symbol for each possible arrangement of the minus signs;
e.g., $X^1-X^2+X^3$, 
$X^1-X^2-X^3$, $X^1+X^2+X^3$
and so on. Let $\Psi^{\pm}$ be the L\'evy exponent for $X^{\pm}_t$. 
Since
$-X^j$ has L\'evy exponent $\overline{\Psi_j},$
$\mathcal{E}_{\Psi^{\pm}}(\mu)=\mathcal{E}_{\Psi}(\mu)$ for all
$X^{\pm}_t$ and
$$\sum\mbox{Re}\left(\int_{\bbR^N_+}
e^{-\sum^N_{j=1}s_j-s\cdot\Psi^{\pm}(\xi)}ds\right)
=2^{N-1}\prod_{j=1}^N\mbox{Re}\left(\frac{1}{1+\Psi_j(\xi)}\right)>0$$
where the first summation $\sum$ is taken over the collection of all the 
$X^{\pm}_t$. On the other hand,
$$Q_{\mu}(\xi)=
\int_{\bbR^N_+}\int_{\bbR^N_+}
e^{-\sum^N_{j=1}|t_j-s_j|\Psi_j(\mbox{sgn}(t_j-s_j)\xi)}\mu(ds)\mu(dt)$$
remains unchanged for all $X^{\pm}_t$ as long as $\mu$ is an $N-$fold 
product measure on $\bbR^N_+.$
Proposition 10.3 of [1] and Theorem 2.1 of [1] together state that
for any additive L\'evy process $X$,
$$k_1\left(\int_{\bbR^d}Q_{\lambda^r}(\xi)d\xi\right)^{-1}\le
E\{\lambda_d(X([0,r]^N))\}\le k_2
\left(\int_{\bbR^d}Q_{\lambda^r}(\xi)d\xi\right)^{-1},$$
where $\lambda^r$ is the restriction of the Lebesgue measure $\lambda_N$ in $\bbR^N$ to
$[0,r]^N$ and $k_1,~k_2\in(0,\infty)$ are two constants depending only
on $r,~N,~d,~\pi.$ Note that $\lambda^r$ is an $N-$fold 
product measure on $\bbR^N_+.$
Thus, 
there exists a constant $c_2\in(0,\infty)$ depending only on
$N$ and $r$ such that 
$$E\{\lambda_d(X([0,r]^N))\}\le c_2
E\{\lambda_d(X^{\pm}([0,r]^N))\}$$ 
for all $X^{\pm}_t$.
Since 
$|1+z|=|1+\bar{z}|$ where $z$ is a complex number,
$\int_{\bbR^d}\prod_{j=1}^N|1+\Psi^{\pm}_j(\xi)|^{-1}d\xi<+\infty$ as well.
Therefore, by (1.4),
\begin{eqnarray*}
&&2^{N-1}\sqrt{c_2}\sqrt{\frac{\mathcal{E}_{\Psi}(\delta_0)}{E\{\lambda_d(X([0,r]^N))\}}}\\
&&\ge\sum\sqrt{\frac{\mathcal{E}_{\Psi^{\pm}}(\delta_0)}{E\{\lambda_d(X^{\pm}([0,r]^N))\}}}\\
&&\ge\sqrt{c_1}\sum\left|\int_{\bbR^d}
\mbox{Re}\left(\int_{\bbR^N_+}
e^{-\sum^N_{j=1}s_j-s\cdot\Psi^{\pm}(\xi)}ds\right)d\xi\right|\\
&&\ge\sqrt{c_1}\left|\sum\int_{\bbR^d}
\mbox{Re}\left(\int_{\bbR^N_+}
e^{-\sum^N_{j=1}s_j-s\cdot\Psi^{\pm}(\xi)}ds\right)d\xi\right|\\
&&=2^{N-1}\sqrt{c_1}\int_{\bbR^d}\prod_{j=1}^N\mbox{Re}\left(\frac{1}{1+\Psi_j(\xi)}\right)d\xi\\
&&=2^{N-1}\sqrt{c_1}(2\pi)^d\mathcal{E}_{\Psi}(\delta_0).
\end{eqnarray*}
(1.3) follows, so does the theorem.~~~~
~~~~~~~~~~~~~~~~~~~~~~~~~~~~~~~~~~~~~~~~~~~~~~~~
~~~~~~~~~~~~~~~~~~~~~~~~~~~~$\Box$

\vskip .25in

\centerline{\textsc{2. Applications}}
\vskip .15in

\noindent
{\it 2.1 The Range of An Additive L\'evy Process }

As the first application, we use Theorem 1.1 to compute 
$\dim_HX(\bbR^N_+).$ Here, $\dim_H$ denotes the Hausdorff dimension.
To begin, we introduce
the standard 
$d$-parameter additive $\alpha$-stable
L\'evy process in $\bbR^d$ for $\alpha\in(0,1):$
$$S^{\alpha}_t=S^1_{t_1}+S^2_{t_2}+\cdots+S^d_{t_d},$$  
that is,
the $S^j$ are independent standard $\alpha$-stable
L\'evy processes in $\bbR^d$
with the common L\'evy exponent
$|\xi|^{\alpha}.$

\vskip 0.15in
\noindent
{\bf Theorem 2.1}~~
{\it Let $X$ be any $N$-parameter additive L\'evy process in $\bbR^d$ with
L\'evy exponent $(\Psi_1,\cdots,\Psi_N)$. Then}
$$\dim_HX(\bbR^N_+)=\sup\left\{\beta\in(0,d):\int_{\bbR^d}|\xi|^{\beta-d}
\prod_{j=1}^N\mbox{Re}\left(\frac{1}{1+\Psi_j(\xi)}\right)
d\xi<\infty\right\}~~~a.s.\eqno(2.1)$$

\vskip 0.15in
\noindent
{\bf Proof}~~
Let $\mathcal{C}_{\beta}$
denote the Riesz capacity. By
Theorem 7.2 of [1], 
for all $\beta\in(0,d)$ and $S^{1-\beta/d}$ independent of $X,$
$$E\mathcal{C}_{\beta}(X(\bbR^N_+))>0\Longleftrightarrow
E\{\lambda_d(S^{1-\beta/d}(\bbR^d_+)+X(\bbR^N_+))\}>0.\eqno(2.2)$$ 
Note that $S^{1-\beta/d}+X$ is a $(d+N,~d)-$additive L\'evy process. 
Thus, by 
Theorem 1.1 and the fact that $\beta<d$ and $\mbox{Re}\left(\frac{1}{1+\Psi_j(\xi)}\right)\in(0,1]$,
we have for all $\beta\in(0,d),$
$$E\mathcal{C}_{\beta}(X(\bbR^N_+))>0\Longleftrightarrow
\int_{\bbR^d}|\xi|^{\beta-d}
\prod_{j=1}^N\mbox{Re}\left(\frac{1}{1+\Psi_j(\xi)}\right)
d\xi<\infty.\eqno(2.3)$$
Thanks to the Frostman theorem,
it remains to show that
$\mathcal{C}_{\beta}(X(\bbR^N_+))>0$ is a trivial event. 
Let $\mathcal{E}_{\beta}$ denote the Riesz energy.
By Plancherel's theorem, given any $\beta\in(0,d),$
there is a constant
$c_{d,\beta}\in(0,\infty)$ such that
$$\mathcal{E}_{\beta}(\nu)=c_{d,\beta}\int_{\bbR^d}|
\hat{\nu}(\xi)|^2|\xi|^{\beta-d}d\xi\eqno(2.4)$$
holds for all probability measures $\nu$
in $\bbR^d.$ Consider the
$1$-killing occupation measure
$$O(A)=\int_{\bbR^N_+}1(X_t\in A)e^{-\sum_{j=1}^Nt_j}dt,~~~A\subset\bbR^d.$$
Clearly, $O$ is a probability measure supported on
$X(\bbR^N_+).$ It is easy to verify that
$$E|\widehat{O}(\xi)|^2=
\prod_{j=1}^N\mbox{Re}\left(\frac{1}{1+\Psi_j(\xi)}\right).$$
It follows from (2.4) that
$$E\mathcal{E}_{\beta}(O)=
c_{d,\beta}\int_{\bbR^d}|\xi|^{\beta-d}
\prod_{j=1}^N\mbox{Re}\left(\frac{1}{1+\Psi_j(\xi)}\right)
d\xi<\infty$$ 
when $E\mathcal{C}_{\beta}(X(\bbR^N_+))>0.$
Therefore,
$\mathcal{E}_{\beta}(O)<\infty$ a.s. Hence,  
$\mathcal{C}_{\beta}(X(\bbR^N_+))>0$ a.s.
~~~~~~~~~~~~~~~~~~~~~~~$\Box$

\vskip .15in
\noindent
{\it 2.2 The Set of $k$-Multiple Points}

First,
we mention a $q$-potential density criterion:
Let $X$ be an additive L\'evy process and assume that
$X$ has an a.e. positive $q$-potential density on $\bbR^d$ for some $q\ge0.$ Then
for all Borel sets $F\subset\bbR^d,$
$$P\left\{
F\bigcap X((0,\infty)^N)\neq\emptyset\right\}>0\Longleftrightarrow
E\left\{\lambda_d(F-X((0,\infty)^N))\right\}>0.\eqno(2.5)$$
The argument is elementary
but crucially hinges on the
property: 
$X_{b+t}-X_b,~t\in\bbR^N_+$ (independent of $X_b$) can be replaced by $X$
for all $b\in\bbR^N_+$; moreover, the second condition ``a.e. positive on $\bbR^d$"
is absolutely necessary for the direction $\Longleftarrow$ in (2.5); 
see for example
Proposition 6.2 of [1].

Let $X^1,\cdots,X^k$ be $k$ independent L\'evy process in $\bbR^d.$
Define
$$Z_t=(X^2_{t_2}-X^1_{t_1},\cdots,X^k_{t_k}-X^{k-1}_{t_{k-1}}),~~~~
t=(t_1,t_2,\cdots,t_k)\in\bbR^k_+.$$
$Z$ is a $k$-parameter additive L\'evy process taking values
in $\bbR^{d(k-1)}.$

\vskip 0.15in
\noindent
{\bf Theorem 2.2}~~
{\it
Let $(X^1;~\Psi_1),~\cdots,~(X^k;~\Psi_k)$ be $k$ independent L\'evy processes in $\bbR^d$ for
$k\ge 2$. 
Assume that $Z$ has an a.e. positive $q$-potential density for some $q\ge0.$ [A special case
is that if 
for each $j=1,\cdots,k,$
$X^j$ has a one-potential density $u_j^1>0,$ $\lambda_d$-a.e., then
$Z$ has an a.e. positive $1$-potential density on $\bbR^{d(k-1)}.$]  
Then}
$$P(\bigcap^k_{j=1}X^j((0,\infty))\neq\emptyset)>0
\Longleftrightarrow
\int_{\bbR^{d(k-1)}}
\prod_{j=1}^k\mbox{Re}\left(\frac{1}{1+\Psi_j(\xi_j-\xi_{j-1})}\right)
d\xi_1\cdots d\xi_{k-1}<\infty\eqno(2.6)$$
{\it with $\xi_0=\xi_k=0.$}

\vskip 0.15in
\noindent
{\bf Proof}~~
For any $\bbR^d$-valued random variable $X$ and $\xi_1,~\xi_2\in\bbR^d,$
$e^{i[(\xi_1,\xi_2)\cdot(X,-X)]}
=e^{i(\xi_1-\xi_2)\cdot X}.$ In particular,
the L\'evy process $(X^j,-X^j)$ has L\'evy exponent
$\Psi_j(\xi_1-\xi_2).$ It follows that
the corresponding integral in (1.1) for $Z$ equals
$$\int_{\bbR^{d(k-1)}}
\prod_{j=1}^k\mbox{Re}\left(\frac{1}{1+\Psi_j(\xi_j-\xi_{j-1})}\right)
d\xi_1\cdots d\xi_{k-1}$$
with $\xi_0=\xi_k=0.$ Clearly,
$$P(\bigcap^k_{j=1}X^j((0,\infty))\neq\emptyset)>0
\Longleftrightarrow
P(0\in Z((0,\infty)^k))>0.$$
Since $Z$ has an a.e. positive $q$-potential density,
by (2.5) 
$$P(0\in Z((0,\infty)^k))>0
\Longleftrightarrow
E\{\lambda_{d(k-1)}(Z((0,\infty)^k))\}>0.$$
(2.6) now follows from Theorem 1.1.
~~~~~~~~~~~~~~~~~~~~~~~~~~~~~~~~~~~~~~~~~~~~~~~~~~~~~~~
~~~~~~~~~~~~~~~~~~~~~~~~~~~~~~~~~~$\Box$

\vskip .15in
For each $\beta\in(0,d)$ and $S^{1-\beta/d}$ independent of $X^1,\cdots,X^k$, define
$$Z^{S,\beta}_t=(X^1_{t_1}-S_{t_0}^{1-\beta/d},X^2_{t_2}-X^1_{t_1},\cdots,X^k_{t_k}-X^{k-1}_{t_{k-1}}),~~
t=(t_0, t_1,t_2,\cdots,t_k)\in\bbR^{d+k}_+,~~t_0\in\bbR^d_+.$$
$Z^{S,\beta}$ is a $k+d$ parameter additive L\'evy process taking values
in $\bbR^{dk}.$

\vskip 0.15in
\noindent
{\bf Theorem 2.3}~~
{\it Let $(X^1;~\Psi_1),~\cdots,~(X^k;~\Psi_k)$ be $k$ independent L\'evy processes in $\bbR^d$ for 
$k\ge 2.$ 
Assume that for each $\beta\in(0,d),$ $Z^{S,\beta}$ has
an a.e. positive $q$-potential density on $\bbR^{dk}$ for some $q\ge0.$ ($q$ might depend on $\beta$.)
[A special case
is that if 
for each $j=1,\cdots,k,$
$X^j$ has a one-potential density $u_j^1>0,$ $\lambda_d$-a.e., then
$Z^{S,\beta}$ has an a.e. positive $1$-potential density on $\bbR^{dk}$ for all $\beta\in(0,d).$]  
If 
$P(\bigcap^k_{j=1}X^j((0,\infty))
\neq\emptyset)>0$, then almost surely 
$\dim_H\bigcap^k_{j=1}X^j((0,\infty))
$ is a constant  
on $\{\bigcap^k_{j=1}X^j((0,\infty))
\neq\emptyset\}$ and} 
$$\dim_H\bigcap^k_{j=1}X^j((0,\infty))
=
\sup\left\{\beta\in(0,d):\int_{\bbR^{dk}}|\sum_{j=1}^k\xi_j|^{\beta-d}
\prod_{j=1}^k\mbox{Re}\left(\frac{1}{1+\Psi_j(\xi_j)}\right)
d\xi_1d\xi_2\cdots d\xi_k<\infty
\right\}.\eqno(2.7)$$

\vskip 0.15in
\noindent
{\bf Proof} ~~
According to the argument, Eq. (4.96)-(4.102), in {\it Proof of Theorem 3.2.} of
Khoshnevisan, Shieh, and Xiao [2], it suffices to show that
for all $\beta\in(0,d)$ and $S^{1-\beta/d}$ independent of
$X^1,\cdots,X^k$, 
$$P\left[\bigcap^k_{j=1}X^j((0,\infty))
\bigcap S^{1-\beta/d}((0,\infty)^d)\neq\emptyset\right]>0\Longleftrightarrow$$
$$\int_{\bbR^{dk}}|\sum_{j=1}^k\xi_j|^{\beta-d}
\prod_{j=1}^k\mbox{Re}\left(\frac{1}{1+\Psi_j(\xi_j)}\right)
d\xi_1d\xi_2\cdots d\xi_k<\infty.\eqno(2.8)$$
Similarly, 
the corresponding integral in (1.1) for $Z^{S,\beta}$ equals
$$\int_{\bbR^{dk}}\frac{1}{(1+|\xi_0|^{1-\beta/d})^d}
\prod_{j=1}^k\mbox{Re}\left(\frac{1}{1+\Psi_j(\xi_j-\xi_{j-1})}\right)
d\xi_0d\xi_1\cdots d\xi_{k-1}$$
with $\xi_k=0.$ Since $Z^{S,\beta}$ has an a.e. positive $q$-potential density,
by (2.5) and Theorem 1.1 
$$P\left[\bigcap^k_{j=1}X^j((0,\infty))
\bigcap S^{1-\beta/d}((0,\infty)^d)\neq\emptyset\right]>0
\Longleftrightarrow
P(0\in Z^{S,\beta}((0,\infty)^{k+d}))>0$$
$$\Longleftrightarrow
E\{\lambda_{dk}(Z^{S,\beta}((0,\infty)^{k+d}))\}>0\Longleftrightarrow$$
$$\int_{\bbR^{dk}}\frac{1}{(1+|\xi_0|^{1-\beta/d})^d}
\prod_{j=1}^k\mbox{Re}\left(\frac{1}{1+\Psi_j(\xi_j-\xi_{j-1})}\right)
d\xi_0d\xi_1\cdots d\xi_{k-1}<\infty$$
with $\xi_k=0.$ Since
$\beta<d$ and 
$\prod_{j=1}^k\mbox{Re}\left(\frac{1}{1+\Psi_j(\xi_j-\xi_{j-1})}\right)\le1,$
$$\int_{\bbR^{dk}}\frac{1}{(1+|\xi_0|^{1-\beta/d})^d}
\prod_{j=1}^k\mbox{Re}\left(\frac{1}{1+\Psi_j(\xi_j-\xi_{j-1})}\right)
d\xi_0d\xi_1\cdots d\xi_{k-1}<\infty$$
$$\Longleftrightarrow
\int_{\bbR^{dk}}|\xi_0|^{\beta-d}
\prod_{j=1}^k\mbox{Re}\left(\frac{1}{1+\Psi_j(\xi_j-\xi_{j-1})}\right)
d\xi_0d\xi_1\cdots d\xi_{k-1}<\infty.$$
Finally,
use the cyclic transformation:
$\xi_j-\xi_{j-1}=\xi'_j,$
$j=1,\cdots,k-1,$ 
$\xi_{k-1}=\xi'_k$ to obtain
$$\int_{\bbR^{dk}}|\xi_0|^{\beta-d}
\prod_{j=1}^k\mbox{Re}\left(\frac{1}{1+\Psi_j(\xi_j-\xi_{j-1})}\right)
d\xi_0d\xi_1\cdots d\xi_{k-1}<\infty$$
$$\Longleftrightarrow
\int_{\bbR^{dk}}|\sum_{j=1}^k\xi'_j|^{\beta-d}
\prod_{j=1}^k\mbox{Re}\left(\frac{1}{1+\Psi_j(\xi'_j)}\right)
d\xi'_1d\xi'_2\cdots d\xi'_k<\infty.\eqno\Box$$

\vskip .15in

Let $X$ be a L\'evy process in $\bbR^d.$
Fix any path $X_t(\omega).$ 
A point $x^{\omega}\in\bbR^d$ is said to be a $k$-multiple point
of $X(\omega)$
if there exist $k$ distinct times $t_1,t_2,\cdots,t_k$ such that
$X_{t_1}(\omega)=X_{t_2}(\omega)=\cdots=X_{t_k}(\omega)=x^{\omega}.$ 
Denote by $E_k^{\omega}$ the set of $k$-multiple points of $X(\omega)$.
It is well known that $E_k$ can be identified with $\bigcap^k_{j=1}X^j((0,\infty))$ 
where the $X^j$ are i.i.d. copies of $X$.
Thus, Theorem 2.2 and Theorem 2.3 imply the next theorem.

\vskip 0.15in
\noindent
{\bf Theorem 2.4}~~ 
{\it Let $(X,~\Psi)$ be any L\'evy process in $\bbR^d.$
Assume that $X$ has a one-potential density $u^1>0,$ $\lambda_d$-a.e.
Let
$E_k$ be the $k$-multiple-point set of $X$. Then} 
$$P(E_k\neq\emptyset)>0\Longleftrightarrow 
\int_{\bbR^{d(k-1)}}
\prod_{j=1}^k\mbox{Re}\left(\frac{1}{1+\Psi(\xi_j-\xi_{j-1})}\right)
d\xi_1\cdots d\xi_{k-1}<\infty\eqno(2.9)$$
{\it with $\xi_0=\xi_k=0.$
If $P(E_k\neq\emptyset)>0,$ then
almost surely $\dim_HE_k$ is a constant   
on $\{E_k\neq\emptyset\}$ and } 
$$\dim_HE_k=
\sup\left\{\beta\in(0,d):\int_{\bbR^{dk}}|\sum_{j=1}^k\xi_j|^{\beta-d}
\prod_{j=1}^k\mbox{Re}\left(\frac{1}{1+\Psi(\xi_j)}\right)
d\xi_1d\xi_2\cdots d\xi_k<\infty
\right\}.\eqno(2.10)$$

\vskip .15in\noindent
{\it 2.3 Intersection of Two Independent Subordinators }

Let $X_t,~t\ge0$ be a process with $X_0=0$, taking values in $\bbR_+.$ First,
we ask this question: What is a condition on $X$ such that for all sets $F\subset(0,\infty),$
$$P(F\bigcap X((0,\infty))\neq\emptyset)>0\Longleftrightarrow
E\{\lambda_1(F-X((0,\infty)))\}>0~~?$$
For subordinators, still the existence and positivity of a $q$-potential density ($q\ge0$)
is the only known useful condition to this question.

Let $\sigma$ be a subordinator. Take an independent copy $\sigma^-$ of $-\sigma.$
We then define a process $\tilde{\sigma}$ on $\bbR$ by
$\tilde{\sigma}_s=\sigma_s$ for $s\ge0$ and $\tilde{\sigma}_s=
\sigma^-_{-s}$ for $s<0.$ Note that $\tilde{\sigma}$ is a process of the property:
$\tilde{\sigma}_{t+b}-\tilde{\sigma}_b,~t\ge0$ (independent of $\tilde{\sigma}_b$)
can be replaced by $\sigma$ for all $b\in\bbR.$

Let $X_t,~t\ge0$ be any process in $\bbR^d$. Then the $q$-potnetial
density is nothing but the density of the expected $q$-occupation measure
with respected to the Lebesgue measure. (When $q=0$, assume that
the expected $0$-occupation measure is finite on the balls.)
Since the reference measure is Lebesgue, one can
easily deduce that if $u$ is a $q$-potential density of $X$,
then $u(-x)$ is a $q$-potential density of $-X$.
Consequently, if we define  
$\widetilde{X}_s=X_s$ for $s\ge0$ and $\widetilde{X}_s=
X^-_{-s}$ for $s<0$ where 
$X^-$ is an independent copy of $-X,$
then $u(x)+u(-x)$ is a $q$-potential density of
$\widetilde{X}.$ Conversely,
if $\widetilde{X}$ has a $q$-potential density, then it
has to be the form $u(x)+u(-x)$, where $u$ is a $q$-potential density
of $X$. If $\sigma$ is a subordinator, after a little thought we can conclude  
that $\tilde{\sigma}$ has an a.e. positive $q$-potential density on $\bbR$ 
if and only if $\sigma$ has an a.e. positive $q$-potential density on $\bbR_+$.

\vskip .15in
\noindent
{\bf Lemma 2.5}~~
{\it If a subordinator $\sigma$ has an a.e. positive $q$-potential density for some $q\ge 0$ on $\bbR_+$,
then for all Borel sets $F\subset(0,\infty)$,}
$$P(F\bigcap \sigma((0,\infty))\neq\emptyset)>0\Longleftrightarrow
E\{\lambda_1(F-\sigma((0,\infty)))\}>0.\eqno(2.11)$$

\vskip .15in
\noindent
{\bf Proof}~~Define $F^*=(-F)\bigcup F.$ Since $F\subset(0,\infty)$,
$(-F)\bigcap F=\emptyset.$ Since $\sigma((0,\infty))\bigcap\tilde{\sigma}((-\infty,0))=\emptyset$ or at most $\{0\},$
by looking at the law of $\sigma^-,$
it is clear that
$$P(F\bigcap \sigma((0,\infty))\neq\emptyset)>0\Longleftrightarrow
P(F^*\bigcap \tilde{\sigma}(\bbR\backslash\{0\})\neq\emptyset)>0.$$
Assume that $E\{\lambda_1(F-\sigma((0,\infty)))\}>0.$ Since
$F\subset F^*,$
$E\{\lambda_1(F^*-\sigma((0,\infty)))\}>0.$ 
From the above discussion, $\tilde{\sigma}$ has an a.e. positive
$q$-potential density. 
Moreover, $\tilde{\sigma}$ is a process of the property:
$\tilde{\sigma}_{t+b}-\tilde{\sigma}_b,~t\ge0$ (independent of $\tilde{\sigma}_b$)
can be replaced by $\sigma$ for all $b\in\bbR.$
It follows from the standard $q$-potential density argument
that 
$P(F^*\bigcap \tilde{\sigma}(\bbR\backslash\{0\})\neq\emptyset)>0.$
The direction $\Longrightarrow$ in (2.11) is elementary since
$\sigma$ has a $q$-potential density.~~~~~~~~~~~~~~
~~~~~~~~~~~~~~~~~~~~~~~~~~~~~~~~~~~~~~~~~~~~~~~~~~$\Box$

\vskip .15in
\noindent
{\bf Theorem 2.6}~~
{\it Let $\sigma^1$ and $\sigma^2$ be two independent subordinators
having the L\'evy exponents $\Psi_1$ and $\Psi_2$, respectively.
Assume that $\sigma^1$ has an a.e. positive $q$-potential
density for some $q\ge0$ on $\bbR_+.$
Then}
$$P[\sigma^1((0,\infty))\bigcap
\sigma^2((0,\infty))\neq\emptyset]>0\Longleftrightarrow
\int_{-\infty}^{\infty}
\mbox{Re}\left(\frac{1}{\Psi_1(x)}\right)
\mbox{Re}\left(\frac{1}{1+\Psi_2(x)}\right)dx<\infty.\eqno(2.12)$$

Note that our result does not require any continuity condition on the $q$-potential density.

\vskip .15in
\noindent
{\bf Proof}~~
By Lemma 2.5 and Theorem 1.1,
$$P[\sigma^1((0,\infty))\bigcap
\sigma^2((0,\infty))\neq\emptyset]>0\Longleftrightarrow
\int_{-\infty}^{\infty}
\mbox{Re}\left(\frac{1}{1+\Psi_1(x)}\right)
\mbox{Re}\left(\frac{1}{1+\Psi_2(x)}\right)dx<\infty.$$
Since $\sigma^1$ is transient,
$\int_{|x|\le1}\mbox{Re}\left(\frac{1}{\Psi_1(x)}\right)dx<\infty.$
The proof is therefore completed.~~~~~~~~~~~~~~~~~~~~~$\Box$

\vskip .15in\noindent
{\it 2.4 A Fourier Integral Problem }

This part of content can be found in Section 6 of [1]. It is
an independent Fourier integral problem. Neither computing
the Hausdorff dimension nor proving the existence of 
$1$-potential density needs the discussion below.
[But this Fourier integral problem might be of novelty to
those who
want to replace the L\'evy exponent
by the $1$-potential density.]   
Let $X$ be an additive L\'evy process. Here is the question.
Suppose that $K:\bbR^d\rightarrow[0,\infty]$ is a symmetric
function with $K(x)<\infty$ for $x\neq0$ that satisfies
$K\in L^1$ and
$\widehat{K}(\xi)=k_1
\prod_{j=1}^N\mbox{Re}\left(\frac{1}{1+\Psi_j(\xi)}\right).$
Under what conditions, can
$$\int\int K(x-y)\mu(dx)\mu(dy)=
k_2\int|\hat{\mu}(\xi)|^2
\prod_{j=1}^N\mbox{Re}\left(\frac{1}{1+\Psi_j(\xi)}\right)d\xi\eqno(2.13)$$
hold for all probability measures $\mu$ in $\bbR^d$?
Here, $k_1,~k_2\in(0,\infty)$ are two constants.
Consider the function $K$ in the following example.
Define $\widetilde{X}^j_{t_j}=-Y^j_{-t_j}$ for $t_j<0$ and $\widetilde{X}^j_{t_j}=X^j_{t_j}$ for $t_j\ge0,$
where $Y^j$ is an independent copy of $X^j$ and
the $Y^j$ are independent of each other and of $X$ as well.
Then
$\widetilde{X}_t=\widetilde{X}^1_{t_1}+\widetilde{X}^2_{t_2}+\cdots+\widetilde{X}^N_{t_N},~
t\in\bbR^N$ is a random field on $\bbR^N.$ Assume that
$\widetilde{X}$ has a $1$-potential density $K$. So, $K\in L^1$ and
a direct check verifies that $K$ is symmetric.
By the definition of $K$,
$\widehat{K}(\xi)=\int_{\bbR^N}
e^{-\sum^N_{j=1}|t_j|}Ee^{i\xi\cdot\widetilde{X}_t}dt.$
Evaluating this integral quadrant by quadrant and using
the identity
$\sum\prod_{j=1}^N\frac{1}{1+z_j^{\pm}}
=2^N\prod_{j=1}^N\mbox{Re}\left(\frac{1}{1+z_j}\right)$ for Re$(z_j)\ge0$
(where $\sum$ is taken over the $2^N$ permutations of conjugate) yield
$\widehat{K}(\xi)=k_1
\prod_{j=1}^N\mbox{Re}\left(\frac{1}{1+\Psi_j(\xi)}\right)>0.$

If $\widehat{K}\in L^1$ (even though this case is less interesting),
on one hand by Fubini,
$$\int|\hat{\mu}(\xi)|^2\widehat{K}(\xi)d\xi
=\int\int\int e^{-i\xi\cdot(x-y)}\widehat{K}(\xi)d\xi\mu(dx)\mu(dy)$$
and on the other hand by inversion 
(assuming the inversion holds everywhere by modification on a null set),
$$\int\int K(x-y)\mu(dx)\mu(dy)=
(2\pi)^{-d}
\int\int\int e^{-i\xi\cdot(x-y)}\widehat{K}(\xi)d\xi\mu(dx)\mu(dy).$$
Thus, (2.13) holds automatically in this case.
If $K$ is continuous at $0$ and $K(0)<\infty$, then $\widehat{K}\in L^1.$
This is a standard fact. Since $K\in L^1$ and
$\widehat{K}>0,$ a bottom line condition needed to prove (2.13) is
that $K$ is continuous at $0$ on $[0,\infty].$
This paper makes no attempt to solve the general case $K(0)=\infty.$

\vskip 0.15in
\noindent
{\bf Remark}~~
Lemma 6.1 of [1] is not valid. (The authors of [1] looked like not having a clear idea how to
prove a result of that sort.) The assumption that
$\mbox{Re}\left(\prod^N_{j=1}\frac{1}{1+\Psi_j(\xi)}\right)>0$
cannot (by any means) justify either equation in (6.4) of [1].
Fortunately, Lemma 6.1 played no role in [1], because
Theorem 7.2 of [1] is an immediate consequence of the
well-known identity (2.4) of the present paper and
Theorem 1.5 of [1]. Nevertheless [1] indeed showed that
the $1$-potential density of an isotropic stable additive
process is comparable to the Riesz kernel at $0$, and therefore
the $1$-potential density is continuous at $0$ on $[0,\infty].$

\vskip .3in

\centerline{REFERENCES}
\vskip .2in

\noindent {[1]} Khoshnevisan, D., Xiao, Y. and Zhong, Y. (2003). 
Measuring the range of an additive 

\indent \indent
L\'evy process. {\it Ann. Probab.} {\bf 31}, 1097-1141.

\noindent {[2]} Khoshnevisan, D., Sheih, N.-R. and Xiao, Y. (2006).
Hausdorff dimension of the contours 

\indent \indent
of symmetric additive proceeses. {\it Probab. Th. Rel. Fields}. (To appear.)

\end{document}